# A RAMSEY PROBLEM RELATED TO BUTTERFLY GRAPH VS. PROPER CONNECTED SUBGRAPHS OF $K_4$

Chula Jayawardene and Lilanthi Samarasekara

*Department of Mathematics, University of Colombo,
Sri Lanka*

e-mails: c_jayawardene@maths.cmb.ac.lk; lilanthi@maths.cmb.ac.lk

A graph on 5 vertices consisting of 2 copies of the cycle graph $C_3$ sharing a common vertex is called the Butterfly graph (*B*). The smallest natural number *s* such that any two-colouring (say red and blue) of the edges of $K_{j \times s}$ has a copy of a red *B* or a blue *G* is called the multipartite Ramsey number of Butterfly graph versus *G*. This number is denoted by $m_j(B,G)$. In this paper we find exact the values for $m_j(B,G)$ when $j \geq 3$ and *G* represents any connected proper subgraph of $K_4$ with at least one edge.

## 1. INTRODUCTION

In this paper we concentrate on simple graphs. Let the complete multipartite graph having *j* uniform sets of size *s* be denoted by $K_{j \times s}$. Given two graphs *G* and *H*, we say that $K_N \rightarrow (G, H)$ if $K_N$ is coloured by two colours red and blue and it contains a copy of *G* (in the first color red) or a copy of *H* (in the second color blue). With respect to this notation, we define the Ramsey number *r(n,m)* as the smallest integer *N* such that $K_N \rightarrow (K_n, K_m)$. As of today, beyond the case *n* = 5, almost nothing significant is known with regard to diagonal classical Ramsey number *r(n,n)* (see [8] for a survey). Burger and Vuuren (see [1]) were honoured for introducing and developing a branch of Ramsey numbers known as size multipartite Ramsey numbers. The size multipartite Ramsey number $m_j(B,G)$, which is a generalization of the much celebrated Ramsey number, is based on exploring the two colourings of multipartite graph $K_{j \times s}$ instead of the complete graph. Formally, we define size multipartite Ramsey number as the smallest natural number *s* such that $K_{j \times s} \rightarrow (K_n, K_m)$.

In the last 14 years, many research papers have been published on the multipartite Ramsey number for different pairs of graphs. [9], has found multipartite Ramsey number for paths versus graph *G* where *G* refers to either a path, a fan, or a windmill. Works of [6,7], focuses on the multipartite Ramsey numbers for graph *G* versus graph *H* where *H* is any isolated vertex free simple graph on four vertices and graph *G* refers to either a $C_3$, a $P_4$ or a $K_4$ - e. This paper presents exact values for $m_j(B,G)$, when $j \geq 3$ where *G* represents a connected proper subgraph of $K_4$ with at least one edge. The details of the results found are summarized in the following table.

# A RAMSEY PROBLEM RELATED TO BUTTERFLY GRAPH VS. PROPER CONNECTED SUBGRAPHS OF $K_4$

| $G =$ | $P_2$ | $P_3$ | $P_4$ | $2K_2$ | $C_3$ | $K_{1,3}$ | $C_4$ | $K_{1,3}+e$ | $B_2$ |
|---|---|---|---|---|---|---|---|---|---|
| $j=3$ | 2 | 2 | 3 | 2 | $\infty$ | 4 | 4 | $\infty$ | $\infty$ |
| $j=4$ | 2 | 2 | 2 | 2 | $\infty$ | 3 | 3 | $\infty$ | $\infty$ |
| $j=5$ | 1 | 1 | 2 | 1 | $\infty$ | 2 | 2 | $\infty$ | $\infty$ |
| $j=6$ | 1 | 1 | 2 | 1 | 2 | 2 | 2 | 2 | 3 |
| $j=7$ | 1 | 1 | 1 | 1 | 2 | 1 | 1 | 2 | 2 |
| $j=8$ | 1 | 1 | 1 | 1 | 2 | 1 | 1 | 2 | 2 |
| $j\geq 9$ | 1 | 1 | 1 | 1 | 1 | 1 | 1 | 1 | 1 |

Table 1.1: Values of $m_j(B,G)$

## 2. NOTATION

Given a graph $G=G(V,E)$ the *order* of the graph is denoted by $|V(G)|$ and the *size* of the graph is denoted by $|E(G)|$. For a vertex $v$ of a graph $G$, the neighbourhood of $v$, denoted by $N(v)$ is defined as the set of vertices adjacent to $v$. Furthermore, the cardinality of this set, denoted $d(v)$, is defined as the degree of $v$. In a Butterfly graph $B$, the vertex of degree 4 is defined as the centre of the Butterfly graph $B$. We say that a graph $G$ is a $k$ regular graph if $d(v) = k$ for all $v \in V(G)$. Let $N_R(v)$ ($N_B(v)$) be the set of vertices adjacent to $v$ in red(blue). Then the cardinality of this set is denoted by $deg_R(v)$ ($deg_B(v)$). Denote the $j$ partite sets of $K_{j\times s}$ by $V_1, V_2, \ldots, V_j$. Let $K_{j\times s} = H_R \oplus H_B$ denote a red and blue coloring of $K_{j\times s}$ where $H_R$ consists of the red graph and where $H_B$ consists of the blue graph, having vertex sets equal to $V(K_{j\times s})$. Suppose that a vertex $u \in V(K_{j\times s})$ of $H_R$ (or $H_B$) belonging to the partite set $V_i$ is such that it is incident to $i_1, i_2, \ldots, i_{j-1}$ vertices of each of the remaining $j$-1 partite sets respectively. Then, we say that vertex $u$ has a $(i_1,i_2,\ldots,i_{j-1})$ red (or blue) split in $H_R$ (or $H_B$) provided that $i_1 \geq i_2 \geq i_3 \geq \ldots \geq i_{j-1}$. Moreover if there exists $k_1,k_2,\ldots,k_{j-1}$ such that $i_1 \geq k_1$, $i_2 \geq k_2$, $i_3 \geq k_3$, $\ldots$ $i_{j-1} \geq k_{j-1}$ and $k_1 \geq k_2 \geq k_3 \geq \ldots \geq k_{j-1}$, then we say that $u$ contains a $(k_1,k_2,\ldots,k_{j-1})$ red (or blue) split in $H_R$ (or $H_B$).

## 3. SIZE RAMSEY NUMBERS FOR $m_j(B,P_2)$ AND $m_j(B,P_3)$

**Theorem 3.1** *If $j \geq 3$, then*

$$m_j(B,P_2) = \begin{cases} 2 & \text{if } j \in \{3,4\} \\ 1 & \text{otherwise} \end{cases}$$

*Proof of Theorem 3.1:* The proof is trivial and is left for the reader.

**Theorem 3.2** *If $j \geq 3$, then*

$$m_j(B,P_2) = \begin{cases} 2 & \text{if } j \in \{3,4\} \\ 1 & \text{otherwise} \end{cases}$$

*Proof of Theorem 3.2:* Consider the red-blue coloring of $K_{3\times 2} = H_R \oplus H_B$ where $H_B$ consists of three independent blue edges $(v_{1,1},v_{2,2})$, $(v_{2,1},v_{3,2})$ and $(v_{1,2},v_{3,1})$. Then $H_R$ will consist of the following diagram. Thus, $K_{3\times 2}$ has neither a blue $P_3$ nor a red $B$. Therefore, $m_3(B,P_3) \geq 3$.

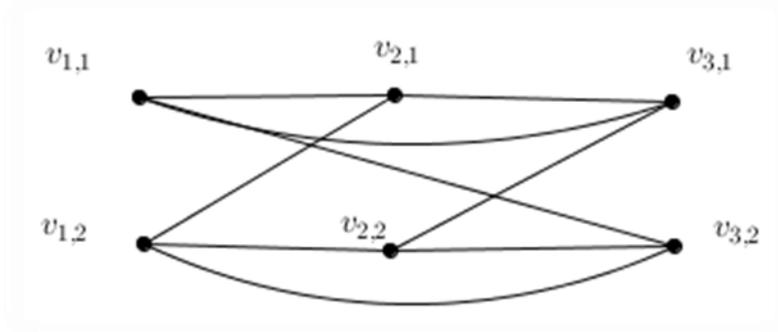

Figure 3.1: The $H_R$ graph

Next to show $m_3(B,P_3) \leq 3$, consider any red/blue coloring given by $K_{3\times 3} = H_R \oplus H_B$ such that $H_R$ contains no red $B$ and $H_B$ contains no blue $P_3$. In order to avoid a blue $P_3$ all vertices must have blue degree at most equal to 1. That is all vertices must have red degree at least equal to 5. Since $K_{3\times 3}$ has odd number of vertices without loss of generality, we may assume that $v_{1,1}$ has red degree 6. However as $deg_R(v_{2,1}) \geq 5$ and $deg_R(v_{2,2}) \geq 5$, there will be a red $2K_2$ induced by $V_2 \cup V_3$. Thus, we will get a red $B$, a contradiction. Therefore, $m_3(B,P_3) = 3$.

As, $r(B,P_3) = 5$, (see [5]) we get $m_4(B,P_3) \geq 2$. Next to show, $m_4(B,P_3) \leq 2$, consider any red/blue coloring given by $K_{4\times 2} = H_R \oplus H_B$, such that $H_R$ contains no red $B$ and $H_B$ contains no blue $P_3$. In order to avoid a blue $P_3$ all vertices must have blue degree at most equal to 1. That is all vertices must have red degree at least equal to 5. Suppose that $v_{1,1}$ is adjacent in red to all vertices of $U=\{v_{2,1}, v_{2,2}, v_{3,1}, v_{3,2}, v_{4,1}\}$. But then in order to avoid a red $B$ induced by $\{v_{1,1}, v_{2,1}, v_{2,2}, v_{3,1}, v_{3,2}, v_{4,1}\}$, $U$ must not contain a red $2K_2$. That is, $U$ must contain a blue $K_{1,2}$, a contradiction. Therefore, $m_4(B,P_3) = 2$.

As, $r(B,P_3) = 5$, we get $m_j(B, P_3) = 1$ for $j \geq 5$.

**Theorem 3.3.** *If $j \geq 3$, then*

$$m_j(B,C_3) = \begin{cases} \infty & j \in \{3,4,5\} \\ 2 & j \in \{6,7,8\} \\ 1 & otherwise \end{cases}$$

*Proof of Theorem 3.3:* $(B,C_3) = \infty$ since $m_j(C_3,C_3) = \infty$ for $j \in \{3,4,5\}$ and $C_3$ is a subgraph of $B$ (See [6]).

Next consider the case $j \in \{6,7,8\}$. First consider a red/blue coloring of $K_{6\times 2}$, given by $K_{6\times 2} = H_R \oplus H_B$, such that $H_R$ contains no red $B$ and $H_B$ contains no blue $C_3$. As $m_6(C_3,C_3) = 1$, the induced subgraph $H_1$ where $V(H_1) = \{v_{i,1} : i \in \{1,2,...,6\}\}$ has a red $C_3$ say $v_{1,1},v_{2,1},v_{3,1},v_{1,1}$. Denote this red $C_3$ by $A_1$. Similarly the induced

# A RAMSEY PROBLEM RELATED TO BUTTERFLY GRAPH VS. PROPER CONNECTED SUBGRAPHS OF $K_4$

subgraph $H_2$ where $V(H_2) = \{v_{3,1}\} \cup \{v_{i,2} : i \in \{1,2,4,5,6\}\}$ has a red $C_3$ (say $A_2$). If $v_{31}$ is a vertex of $A_2$ then $K_{6\times 2}$ has a red $B$, a contradiction. Otherwise, we get the following three cases.

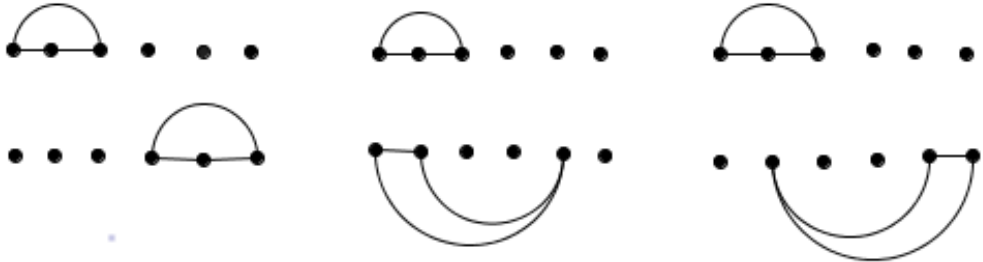

Figure 3.2: The three cases

*Case 1:* None of the vertices of $A_2$ belong to the partite sets $V_1, V_2, V_3$.

*Case 2:* Two of the vertices of $A_2$ belong to two of the partite sets $V_1, V_2, V_3$.

*Case 3:* Only one of the vertices of $A_2$ belong to one of the partite sets $V_1, V_2, V_3$.

In each of these cases first consider the induced subgraph $H$ such that $H$ consists of six vertices $v_1, v_2, ..., v_6$ where no two vertices of $\{v_1, v_2, ..., v_6\}$ belong to the same partite set and $|V(A_i) \cap V(H)| = 1$ for each $i \in \{1,2\}$. Due to the absence of a blue $C_3$ and $m_6(C_3, C_3) = 1$, $H$ has a red $C_3$ (say $A_3$). If $V(A_3) \cap V(A_1) \neq \emptyset$ or $V(A_3) \cap V(A_2) \neq \emptyset$ then $K_{6\times 2}$ has a red $B$, a contradiction.

Otherwise, consider the induced subgraph $H_1$ consisting of the six vertices $u_1, u_2, ..., u_6$ where no two vertices of $\{u_1, u_2, ..., u_6\}$ belong to the same partite set and $|V(A_i) \cap V(H_1)| = 2$ for each $i \in \{1,2,3\}$. $H_1$ has a red $C_3$ due to the absence of a blue $C_3$ and $m_6(C_3, C_3) = 1$. This red $C_3$ along with one of the $A_i$ where $i \in \{1,2,3\}$ forms a red $B$, a contradiction. Therefore, $m_6(B, C_3) \leq 2$.

As, $r(B, C_3) = 9$, (see [5]) we get, $m_8(B, C_3) \geq 2$.

Therefore, $2 \leq m_8(B, C_3) \leq m_7(B, C_3) \leq m_6(B, C_3) \leq 2$, gives us $m_j(B, C_3) = 2$ for $j \in \{6,7,8\}$.

Finally, as $r(B, C_3) = 9$, (see [5]) we get, $m_j(B, C_3) = 1$ if $j \geq 9$. $\square$

**Theorem 3.4.** *If $j \geq 3$, then*

$$m_j(B, P_4) = \begin{cases} 3 & j=3 \\ 2 & j \in \{4,5,6\} \\ 1 & otherwise \end{cases}$$

*Proof of Theorem 3.4:* $m_3(B, P_4) \geq 3$. since $m_3(B, P_3) = 3$ by theorem 2.

Next to show, $m_3(B,P_4) \leq 3$, consider any red/blue coloring of $K_{3\times 3}$ given by $K_{3\times 3} = H_R \oplus H_B$, such that $H_R$ contains no red $B$ and $H_B$ contains no blue $P_4$. By the theorem 2 as $m_3(B,P_3)=3$, we get that there exists a blue $P_3$.

*Case 1:* There exists a blue $P_3$ that lies in three partite sets.

Without loss of generality, assume that this blue $P_3$ comprises of $(v_{1,1},v_{2,1})$ and $(v_{2,1},v_{3,1})$ blue edges. But then in order for $K_{3\times 3}$ not to have a blue $P_4$, $(v_{3,1},v_{1,2})$, $(v_{3,1},v_{1,3})$, $(v_{3,1},v_{2,2})$ and $(v_{3,1},v_{2,3})$ have to be red edges. Next for $W=\{v_{1,2},v_{1,3},v_{2,2},v_{2,3},v_{3,1}\}$ not to induce a red $B$, $W$ will be forced to contain a blue $P_3$, belonging to $V_1$ and $V_2$. Thus, this case leads to the following case 2.

*Case 2:* There exists a blue $P_3$ that lies in two partite sets.

Without loss of generality, assume that this blue $P_3$ comprises of $(v_{1,1},v_{2,1})$ and $(v_{1,1}, v_{2,2})$ blue edges. But then in order for $K_{3\times 3}$ not to have a blue $P_4$, $\{v_{2,1},v_{2,2}\}$ will have to be adjacent to all vertices of $W_1=\{v_{1,2},v_{1,3},v_{3,1},v_{3,2},v_{3,3}\}$ in red. In order for $W_1$ not to induce a blue $P_4$, without loss of generality we may assume that, $(v_{1,2},v_{3,1})$ is a red edge. But then in order to avoid a red $B$ induced by $W_1 \cup \{v_{2,1}\}$, $(v_{1,3}, v_{3,2})$ and $(v_{1,3}, v_{3,3})$ are blue edges. In addition, in order to avoid blue $P_4$, given by $v_{3,2}\ v_{1,3}\ v_{3,3}\ v_{1,2}$ the edge $(v_{3,3}, v_{1,2})$ are a red edge. But then $\{v_{1,2},v_{2,1},v_{2,2},v_{3,1},v_{3,3}\}$ will induce a red $B$ consisting of the two red triangle $v_{1,2},v_{2,2},v_{3,1},v_{1,2}$ and $v_{1,2},v_{2,1},v_{3,3},v_{1,2}$ with $v_{1,2}$ as the centre vertex, a contradiction.

Thus, $m_3(B,P_4) \leq 3$. Therefore, $m_3(B,P_4) = 3$.

As $r(B,P_4) = 7$, (see [5]) we get $m_6(B, P_4) \geq 2$. To show $m_4(B, P_4) \leq 2$, consider $K_{4\times 2}$ with any red/blue coloring. Assume $K_{4\times 2}$ has neither a red $B$ nor a blue $P_4$. Since $m_4(B,P_3) = 2$ and $K_{4\times 2}$ has no red $B$, it has a blue $P_3$.

*Case 1:* There exists a blue $P_3$ that lies in two partite sets

Let the blue $P_3$ be $v_{1,1}$, $v_{2,1}$, $v_{1,2}$. As there is no blue $P_4$ all vertices in $V_3 \cup V_4 \cup \{v_{2,2}\}$ are adjacent in red to both $v_{1,1}$ and $v_{1,2}$. As there is no red $B$ the red graph induced by $H=\{v_{2,2}, v_{3,1}, v_{3,2}, v_{4,1}, v_{4,2}\}$ has no red $2K_2$. Then any connected components in the graph induced by $H$ is equal to a $K_1$, $K_2$, $P_3$, $K_3$, $K_{1,3}$ or $K_{1,4}$. Also the induced red graph of $H$ can contain at most one connected component having one or more red edges. In both these situations the blue graph induced by $H$ has a blue $P_4$, a contradiction.

*Case 2:* There exists a blue $P_3$ that lies in three partite sets.

Let the blue $P_3$ be $v_{1,1}$, $v_{2,1}$, $v_{3,1}$. As there is no blue $P_4$ all vertices in $\{v_{i,2}: i \in \{2,3,4\}\} \cup \{v_{4,1}\}$ are adjacent in red to $v_{1,1}$ and all vertices in $\{v_{i,2} : i \in \{1,2,4\}\} \cup \{v_{4,1}\}$ are adjacent in red to $v_{3,1}$. However, by the elimination of case 1, either $(v_{3,2},v_{4,1})$ or $(v_{3,2},v_{4,2})$ must be red. Without loss of generality assume that $(v_{3,2}, v_{4,1})$ is red. In order to avoid a red $B$ with $v_{1,1}$ as the centre, $(v_{2,2},v_{4,2})$ must be blue. But then, in order to avoid case 1, $(v_{2,2}, v_{4,1})$ is red. Next in order to avoid a red $B$ with $v_{1,1}$ as the centre, $(v_{3,2}, v_{4,2})$ must be blue and in order to avoid a red $B$ with $v_{3,1}$ as the centre, $(v_{1,2}, v_{4,2})$ must be blue. However, as there is no blue $P_4$, $(v_{1,2}, v_{2,2})$, $(v_{2,2}, v_{3,2})$ and $(v_{1,2}, v_{3,2})$ must all be red. This gives us a red $B$ with $v_{2,2}$ as the centre, a contradiction. Therefore, $m_4(B, P_4) \leq 2$.

As $2 \leq m_6(B, P_4) \leq m_5(B, P_4) \leq m_4(B, P_4) \leq 2$ we get $m_j(B, P_4) = 2$ for $j \in \{4,5,6\}$.
Finally, as $r(B, P_4) = 7$ we get, $m_j(B, P_4) = 1$ if $j \geq 7$. □

**Theorem 3.5.** *If $j \geq 3$, then*

$$m_j(B, 2K_2) = \begin{cases} 2 & \text{if } j \in \{3,4\} \\ 1 & \text{otherwise} \end{cases}$$

*Proof of Theorem 3.5.* As $r(B, 2K_2) = 5$, (see [5]), we get $m_4(B, 2K_2) \geq 2$. To show $m_3(B, 2K_2) \leq 2$, consider $K_{3 \times 2}$ with any red/blue coloring. As $m_3(B, P_2) = 2$ and $K_{3 \times 2}$ has no red $B$, it has a blue $P_2$ (say $(v_{1,1}, v_{2,1})$). As there is no blue $2K_2$, $(v_{1,2}, v_{2,2})$, $(v_{1,2}, v_{3,1})$, $(v_{1,2}, v_{3,2})$, $(v_{2,2}, v_{3,1})$ and $(v_{2,2}, v_{3,2})$ and are red edges. Next if either $(v_{1,1}, v_{2,2})$ or $(v_{1,1}, v_{3,1})$ is blue we would get that $V_2 \cup V_3 \cup \{v_{1,2}\}$ will induce a red $B$. Therefore, both $(v_{1,1}, v_{2,2})$ and $(v_{1,1}, v_{3,1})$ will have to be red. However in this case too, we will get that $V_1 \cup \{v_{2,2}\} \cup V_3$ will induce a red $B$, a contradiction. Hence, $m_3(B, 2K_2) \leq 2$. Therefore, $2 \leq m_4(B, 2K_2) \leq m_3(B, 2K_2) \leq 2$, gives us $m_3(B, 2K_2) = 2$ and $m_4(B, 2K_2) = 2$.

Finally, as $r(B, 2K_2) = 5$, we get, $m_j(B, 2K_2) = 1$ if $j \geq 5$. □

**Theorem 3.6.** *If $j \geq 3$, then*

$$m_j(B, K_{1,3}) = \begin{cases} 4 & j = 3 \\ 3 & j = 4 \\ 2 & j \in \{5, 6\} \\ 1 & \text{otherwise} \end{cases}$$

*Proof of Theorem 3.6.* For the case $j = 3$, consider the red-blue coloring of $K_{3 \times 3} = H_R \oplus H_B$ where $H_B$ is the blue cycle $v_{1,1} v_{3,3} v_{1,2} v_{2,1} v_{1,3} v_{2,2} v_{3,1} v_{2,3} v_{3,2} v_{1,1}$ (see Theorem 9: figure 6). Then for any vertex $v \in V(K_{3 \times 3})$, the red induced subgraph of $V(\{v\} \cup N_R(v))$ will be isomorphic to a $K_{1,4} + e$. Therefore, $K_{3 \times 3}$ has neither a blue $K_{1,3}$ nor a red $B$. Therefore, $m_3(B, K_{1,3}) \geq 4$. Consider any red-blue coloring of $K_{3 \times 4} = H_R \oplus H_B$ such that $H_R$ contains no red $B$ and $H_B$ contains no blue $K_{1,3}$. In order to avoid a blue $K_{1,3}$ for any vertex $v \in V(K_{3 \times 4})$, $\deg_B(v) \leq 2$ and $\deg_R(v) \geq 6$. Thus, without loss of generality assume that $v_{1,1}$ is adjacent in red to

$v_{2,1},v_{2,2},v_{2,3}$, $v_{3,1},v_{3,2}$. But then as $deg_R(v_{3,1}) \geq 6$ and $deg_R(v_{3,2}) \geq 6$, we get $v_{3,1}$ and $v_{3,2}$ are adjacent in red to at least one vertex of $\{v_{2,1},v_{2,2},v_{2,3}\}$. Therefore, we get the red edges $(v_{3,1},x)$ and $(v_{3,2},y)$ for some $x, y \in \{v_{2,1},v_{2,2},v_{2,3}\}$ such that $x \neq y$ or else $(v_{3,1},x)$ and $(v_{3,2},x)$ for some $x \in \{v_{2,1},v_{2,2},v_{2,3}\}$. In the first possibility, we get a red $B$ induced by $\{v_{1,1},v_{3,1},x,v_{3,2},y\}$, a contradiction. In the second possibility, $\{v_{2,1},v_{2,2},v_{2,3}\} \cap \{x\}^c$ will be adjacent in blue to all vertices of $\{v_{3,1},v_{3,2}\}$ in blue as we have already eliminated the first possibility. But then as $v_{3,2}$ is adjacent to all vertices of $V_1$ and as $deg_R(x) \geq 6$, we get that $v_{3,2}$ and $x$ will have a common red neighbour say $y$ in $V_1$ distinct from $v_{1,1}$. This will give us a red $B$ induced by $\{v_{1,1},v_{3,1},x,v_{3,2},y\}$ with centre $x$, a contradiction. Therefore, $m_3(B,K_{1,3}) = 4$.

Now let us deal with the case $j = 4$. Consider the red/blue coloring of $K_{4\times 2} = H_R \oplus H_B$ where $H_B$ is a blue $2C_4$ such that $v_{1,1}v_{2,1}v_{1,2}v_2$, $v_{1,1}$ and $v_{3,1}v_{4,1}v_{3,2}v_{4,2}v_{3,1}$ represents blue $C_4$'s. Then as $H_R$ will consist of a red $K_{4,4}$ and hence $K_{4\times 2}$ will neither have a blue $K_{1,3}$ nor a red $B$. Therefore, $m_4(B,K_{1,3}) \geq 3$. Next to show, $m_4(B,K_{1,3}) \geq 3$, consider any red-blue coloring of $K_{4\times 3} = H_R \oplus H_B$ such that $H_R$ contains no red $B$ and $H_B$ contains no blue $K_{1,3}$. In order to avoid a blue $K_{1,3}$ for any vertex $v \in V(K_{4\times 3})$, $deg_B(v) \leq 2$ and $deg_R(v) \geq 7$. Thus, $v_{1,1}$ contains a red (2,2,1) red split and thus without loss of generality we may assume that, $v_{1,1}$ is adjacent in red to $\{v_{2,1},v_{2,2},v_{3,1},v_{3,2},v_{4,1}\}$. As, $v_{4,1}$ contains a red (2,2,1) red split and thus without loss of generality we may assume that, $(v_{3,1},v_{4,1})$ is a red edge. Next to avoid a red $B$, $(v_{2,1},v_{3,2})$ and $(v_{2,2},v_{3,2})$ are blue edges. As $deg_R(v_{3,2}) \geq 7$, we get $(v_{3,2},v_{4,1})$ is a red edge. But then in order to avoid a red $B$, we get $(v_{2,1}, v_{3,1})$ and $(v_{2,2},v_{3,1})$ are blue edges. But $deg_R(v_{1,1}) \geq 7$ and $v_{3,1}$ is adjacent in red to all vertices of $V_1 \cup \{v_{2,3}\} \cup V_4$, we also get $v_{1,1}$ and $v_{3,1}$ will have a common red neighbour distinct from $v_{4,1}$ say $y$. But this will force $\{v_{1,1},y,v_{3,1},v_{3,2},v_{4,1}\}$ to induce a red $B$ with centre, $v_{1,1}$ a contradiction. Hence, $m_4(B,K_{1,3}) \geq 3$. Therefore, $m_4(B,K_{1,3}) = 3$, as required.

Now let us deal with the case $j \in \{5,6\}$. As $C_3$ is a subgraph of $B$, $m_6(B,K_{1,3}) \geq m_6(C_3,K_{1,3})$. Since $m_6(C_3,K_{1,3}) = 2$ ( See [6]), $m_6(,K_{1,3}) \geq 2$. Next consider the red/blue coloring of $K_{5\times 2} = H_R \oplus H_B$. Assume $K_{5\times 2}$ has no red $B$ and has no blue $K_{1,3}$. As $m_5(C_3,K_{1,3}) = 2$, $K_{5\times 2}$ has a red $C_3$ (say $v_{1,1}v_{2,1}v_{3,1}v_{1,1}$). Also, for any vertex $v$, $deg_B(v) \leq 2$ and hence $deg_R(v) \geq 6$. Therefore, $v_{1,1}$ is adjacent to at least 4 vertices in $X = \{v_{4,1},v_{5,1},v_{2,2},v_{3,2},v_{4,2},v_{5,2}\}$. This gives the following 3 situations.

*Situation 1:* The 4 vertices are in 4 partite sets.

*Situation 2:* The 4 vertices are in 3 partite sets. Note there are two possibilities under this situation, namely, situation 2(a) and situation 2(b) as illustrated in the following diagram.

*Situation 3:* The 4 vertices are in 2 partite sets.

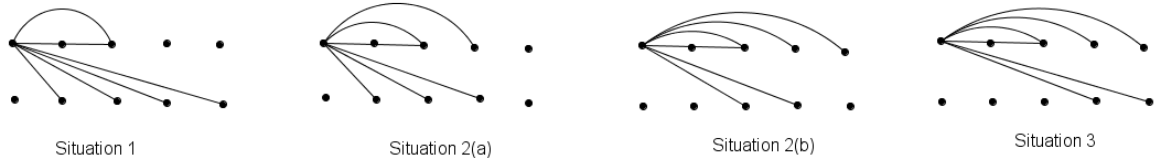

| Situation 1 | Situation 2(a) | Situation 2(b) | Situation 3 |

Figure 3.3: The three possible situations: situations 1 to 3 are from left to right

As $deg_B(v) \leq 2$ for any vertex $v$, one of the edges in $\{(v_{2,2},v_{i,2}) : i \in \{3,4,5\}\}$ is red in situation 1 and one of the edges in $\{(v_{2,2},v_{3,2}),(v_{2,2},v_{4,1}),(v_{2,2},v_{4,2})\}$ is red in situation 2(a) and one of the edges in $\{(v_{3,2},v_{4,1}),(v_{3,2},v_{4,2}),(v_{3,2},v_{5,1})\}$ is red in situation 2(b). This result in $K_{5\times 2}$ having a red $B$, a contradiction.

Consider when $K_{5\times 2}$ is under situation 3. If one of the edges in $\{(v_{4,1},v_{5,1}),(v_{4,1},v_{5,2}),(v_{4,2},v_{5,1}),(v_{4,2},v_{5,2})\}$ is red then $K_{5\times 2}$ has a red $B$, a contradiction. Therefore, all edges in $\{(v_{4,1},v_{5,1}),(v_{4,1},v_{5,2}),(v_{4,2},v_{5,1}),(v_{4,2},v_{5,2})\}$ are blue. But then the edges $(v_{4,1},v_{3,1})$ and $(v_{4,2},v_{2,1})$ are red, due to $deg_B(v) \leq 2$ for any vertex $v$. Then $K_{5\times 2}$ has a red $B$ induced by $\{v_{1,1},v_{2,1},v_{3,1},v_{4,1},v_{4,2}\}$ with centre $v_{1,1}$, a contradiction. Therefore, $m_5(B,K_{1,3}) \leq 2$. Then, as $2 \geq m_5(B,K_{1,3}) \geq m_6(B,K_{1,3}) \geq 2$, we get $m_j(B,K_{1,3}) = 2$ for $j \in \{5,6\}$.

Finally, as $r(B,K_{1,3}) = 7$, we get, $m_j(B,K_{1,3}) = 1$ if $j \geq 7$. □

**Theorem 3.7.** *If $j \geq 3$, then*

$$m_j(B, K_{1,3}+x) = \begin{cases} \infty & j \in \{3,4,5\} \\ 2 & j \in \{6,7,8\} \\ 1 & otherwise \end{cases}$$

*Proof of Theorem 3.7.* $m_j(B,K_{1,3}+x) = \infty$ since $m_j(C_3,C_3) = \infty$ for $j \in \{3,4,5\}$ and $C_3$ is a subgraph of $K_{1,3}+x$ (See [6]). Consider any red-blue colouring of $K_{6\times 2}$. Assume $K_{6\times 2}$ has no red $B$. Then the graph has a blue $C_3$ (say $H$) where $V(H) = \{v_{1,1},v_{2,1},v_{3,1}\}$ as $m_6(B,C_3) = 2$. If one of the vertices in $\{v_{1,1},v_{2,1},v_{3,1}\}$ is incident to a blue edge that is not in $E(H)$ then $K_{6\times 2}$ has a blue $K_{1,3}+x$. Otherwise, all edges in the sets $\{(v_{3,1},v_{1,2}),(v_{3,1},v_{2,2})\}$, $\{(v_{3,1},v_{i,1}) : i \in \{4,5,6\})\}$ and $\{(v_{3,1},v_{i,2}) : i \in \{4,5,6\}\}$ are red. Let $W = \{v_{1,2},v_{2,2}\} \cup \{v_{i,1}: i \in \{4,5,6\})\} \cup \{v_{i,2} : i \in \{4,5,6\})\}$.

Then due to the absence of a red $2K_2$ in $W$, the red graph induced by $W$ can have only one red component of size greater than or equal to one, Also this component will be equal to a red $C_3$ or a $K_{1,n}$ where $1 \leq n \leq 7$. But in all such situations the blue graph induced by $W$ will contain a blue $K_{1,3} + x$. Therefore, $m_6(B,K_{1,3}+x) \leq 2$. As we know that $m_6(B,K_{1,3}+x) \geq m_7(B,K_{1,3}+x) \geq m_8(B,K_{1,3}+x) \geq m_8(B,C_3) = 2$ we get, $m_j(B,K_{1,3}+x) = 2$ for $j \in \{6,7,8\}$.

Finally, as $r(B,C_3) = 9$ (see [5]) we get, $m_j(B,C_3) = 1$ if $j \geq 9$. □

**Theorem 3.8.** *If $j \geq 3$, then*

$$m_j(B, B_2) = \begin{cases} \infty & j \in \{3,4,5\} \\ 3 & j = 6 \\ 2 & j \in \{7,8\} \\ 1 & \text{otherwise} \end{cases}$$

*Proof of Theorem 3.8.* As $m_j(B, C_3) = \infty$ for $j \in \{3,4,5\}$ and $m_j(B, B_2) \geq m_j(B, C_3)$ we get that $m_j(B, B_2) = \infty$ for $j \in \{3,4,5\}$.

In the case $j = 6$, consider the red-blue coloring of $K_{6 \times 2} = H_R \oplus H_B$ where $H_R$ and $H_B$ are the red and blue graphs given in figure 4 and figure 5 respectively.

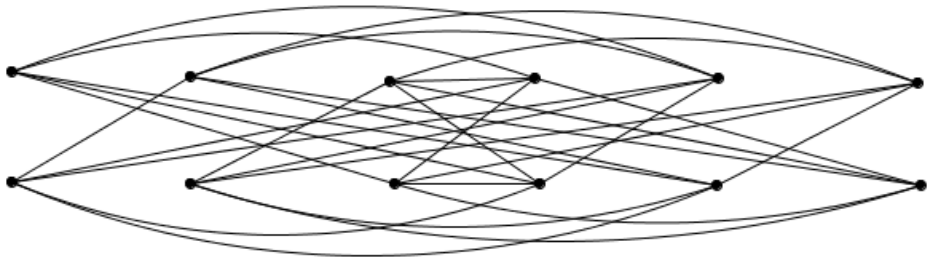

Figure 3.4: $H_R$ graph for $m_6(B, B_2) \geq 3$

From the figure 4, it is evident that $H_R$ contains only 3 disjoint red $B_2$ graphs and thus does not contain a red $B$.

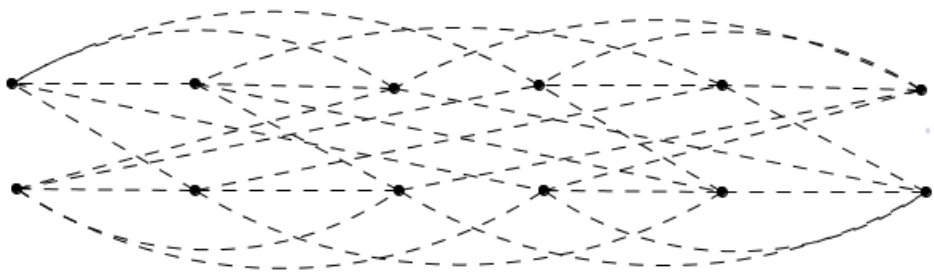

Figure 3.5: $H_B$ graph for $m_6(B, B_2) \geq 3$

From the figure 5, it is evident that $H_B$ contains only 4 disjoint blue triangles and thus does not contain a blue $B_2$. Therefore, we get $m_6(B, B_2) \geq 3$.

# A RAMSEY PROBLEM RELATED TO BUTTERFLY GRAPH VS. PROPER CONNECTED SUBGRAPHS OF $K_4$

Next, consider any red/blue coloring of $K_{6\times 3}= H_R \oplus H_B$ such that $H_R$ contains no red $B$ and $H_B$ contains no blue $B_2$. Then, the graph has a blue $C_3$ (say $v_{1,1}v_{2,1}v_{3,1}v_{1,1}$) as $m_6(B,C_3) = 2$.

If there is a vertex in $Y = \cup^3_{i=1}\{v_{i,l}: l \in \{2,3\}\}$ such that it is adjacent in blue to 2 vertices in $\{v_{1,1},v_{2,1},v_{3,1}\}$, then $K_{6\times 3}$ has a blue $B_2$, a contradiction. Therefore, by Pigeon Hole principle we may assume that one of the vertices of $\{v_{1,1},v_{2,1},v_{3,1}\}$ is adjacent in red to at least two vertices of $Y$ (which may or may belong to one or two partite sets). Next, if there is a vertex in $X = \cup^6_{i=4}\{v_{i,l}: l \in \{1,2,3\}\}$ such that it is adjacent in blue to 2 vertices in $\{v_{1,1},v_{2,1},v_{3,1}\}$. Then $K_{6\times 3}$ has a blue $B_2$, a contradiction. Therefore, all vertices in $X$ are adjacent in red to at least 2 vertices in $\{v_{11},v_{21},v_{31}\}$. Then, there are at least 18 red edges from $X$ to $v_{1,1}v_{2,1}v_{3,1}v_{1,1}$. By Pigeon Hole principal, there exist a vertex in $\{v_{1,1},v_{2,1},v_{3,1}\}$ adjacent to 6 of the red edges from $X$.

*Remark:* If a vertex of $\{v_{1,1},v_{2,1},v_{3,1}\}$ is adjacent to at least 6 vertices in $X \cup Y$ such that 6 vertices belong to 3 partite sets where each partite set has 2 vertices of the 6 vertices each or else is adjacent to 6 vertices in $X \cup Y$ such that 6 vertices belong to 4 partite sets where each of the four partite sets have at least one vertex of the 6 vertices then there exists a blue $B_2$ or red $B$.

By counting at least two vertices of $\{v_{1,1},v_{2,1},v_{3,1}\}$ must be adjacent to at least 6 vertices of $X$ each, as no vertex of $\{v_{1,1},v_{2,1},v_{3,1}\}$ can be adjacent to at 8 or more vertices by the above remark. Without loss of generality assume that these two vertices are $v_{1,1}$ and $v_{2,1}$. By the above remark, $v_{1,1}$ and $v_{2,1}$ must contain a red (3,3,0) or a red (3,2,1) split in $X$.

Without loss of generality if $v_{1,1}$ has a red (3,2,1) split in $X$, by the above remark, then $v_{1,1}$ cannot be adjacent in red to any vertex of $Y$. But then by counting $v_{2,1}$ must be adjacent in red to at least 2 vertices in $Y$. But by remark this will result in a blue $B_2$ or red $B$. Therefore, we are left with the option both $v_{1,1}$ and $v_{2,1}$ must contain a red (3,3,0) split in $X$ and more precisely, equal to a red (3,3,0) split in $X$. However, as each of the vertices of $X$ must be adjacent in red to two vertices of $\{v_{1,1},v_{2,1},v_{3,1}\}$, $v_{3,1}$ will contain a red (3,3,0) split in $X$. Therefore, each of the vertices $v_{1,1},v_{2,1},v_{3,1}$, will contain a red (3,3,0) split in $X$. However, by Pigeon Hole principle as one vertex of $\{v_{1,1},v_{2,1},v_{3,1}\}$ is adjacent in red to at least two vertices of $Y$, the above remark will give a red $B$ or blue $B_2$, a contradiction. Therefore, we get that, $m_6(B,B_2) \leq 3$. That is, $m_3(B,B_2) = 3$, as required.

In the case $j = \{7,8\}$, first consider any red-blue colouring of $K_{7\times 2}= H_R \oplus H_B$ such that $H_R$ contains no red $B$ and $H_B$ contains no blue $B_2$. Then, the graph has a blue $C_3$ (say $v_{1,1}v_{2,1}v_{3,1}v_{1,1}$) as $m_7(B,C_3) = 2$. If there is a vertex in $X = \cup^7_{i=4}\{v_{i,l}: l \in \{1,2\}\}$ such that it is adjacent in blue to 2 vertices in $\{v_{1,1},v_{2,1},v_{3,1}\}$, then $K_{7\times 2}$ has a blue $B_2$, a contradiction. Therefore, all vertices in $X$ are adjacent in red to at least 2 vertices in $\{v_{1,1},v_{2,1},v_{3,1}\}$. Then, there are at least 16 red edges from $X$ to $v_{1,1}v_{2,1}v_{3,1}v_{1,1}$. By Pigeon Hole principal, there exists a vertex in $\{v_{1,1}, v_{2,1}, v_{3,1}\}$ (say $v_{3,1}$) adjacent to 6 of the red edges from $X$. This gives rise to the two situations,

*Situation 1:* The 6 vertices are in 4 partite sets.

This results in a containing a red (2,2,1,1) split in $X$. Since $N_R(v_{3,1})$ cannot contain a red $2K_2$ Then any connected components in the graph induced by $N_R(v_{3,1})$ is equal to a $K_1$, $K_2$, $P_3$, $K_3$, $K_{1,3}$, $K_{1,4}$ or $K_{1,5}$. Also the induced red graph of $H$ can contain at most one connected component having one or more red edges. This will result in a red $B$ or a blue $B_2$, a contradiction.

*Situation 2:* The 6 vertices are in 3 partite sets.

This results in a containing a red (2,2,2,0) split in $X$. Since $N_R(v_{3,1})$ cannot contain a red $2K_2$ Then any connected components in the graph induced by $N_R(v_{3,1})$ is equal to a $K_1$, $K_2$, $P_3$, $K_3$ $K_{1,3}$ or $K_{1,4}$. Also the induced red graph of $H$ can contain at most one connected component having one or more red edges. This will result in a red $B$ or a blue $B_2$, a contradiction.

Therefore, $m_7(B,B_2) \leq 2$. As $m_7(B,B_2) \geq m_8(B,B_2) \geq m_8(B,C_3) = 2$, for $j \in \{7,8\}$ we get $m_j(B,B_2) = 2$. Finally, as $r(B,B_2)=9$ (see [5]) we get, $m_j(B,B_2) = 1$ if $j \geq 9$. □

**Theorem 3.9.** *If $j \geq 3$, then*

$$m_j(B,C_4) = \begin{cases} 4 & j=3 \\ 3 & j=4 \\ 2 & j \in \{5,6\} \\ 1 & \text{otherwise} \end{cases}$$

*Proof of Theorem 3.9.* If $j = 3$, Consider the red-blue coloring of $K_{3\times 3} = H_R \oplus H_B$ where $H_B$ is the blue cycle $v_{1,1}v_{3,3}v_{1,2}v_{2,1}v_{1,3}v_{2,2}v_{3,1}v_{2,3}v_{3,2}v_{1,1}$ as illustrated in the following diagram. Then for any vertex $v \in V(K_{3\times 3})$, the red induced subgraph of $V(\{v\} \cup N_R(v))$ will be isomorphic to a $K_{1,4}+e$. Thus, $K_{3\times 3}$ has neither a blue $C_4$ nor a red $B$. Therefore, we get $m_3(B,C_4) \geq 4$.

Next to show, $m_3(B,C_4) \leq 4$, consider any red/blue coloring given by $G = K_{3\times 4} = H_R \oplus H_B$, such that $H_R$ contains no red $B$ and $H_B$ contains no blue $C_4$.

# A RAMSEY PROBLEM RELATED TO BUTTERFLY GRAPH VS. PROPER CONNECTED SUBGRAPHS OF $K_4$

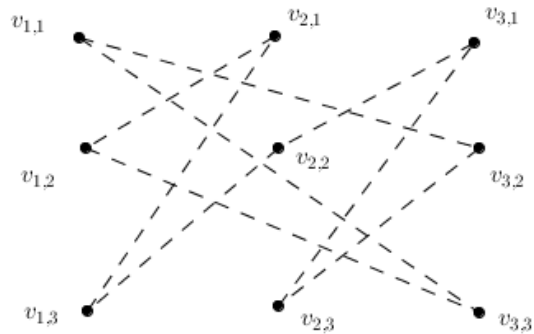

Figure 3.6: The $H_B$ graph used to show $m_3(B,C_4) \geq 4$.

*Claim 1: a)* No vertex of $G$ can have red $(x,y)$ split where $x,y \geq 3$. Therefore, $\delta_B(H_B) \geq 2$.

*b)* No vertex of $G$ can contain a red $(3,2)$ split. Therefore, $\delta_B(H_B) \geq 3$ and if there is a vertex $v$ such that $deg_B(v) = 3$ then it contains a blue $(3,0)$ split.

*Proof of claim 1(a).* Suppose that Claim 1(a) is false. In particular, suppose that $v$ contains a $(3,3)$ split then as $H_R$ has no red $B$, $N_R(v)$ cannot have a red $2K_2$. This would force $N_R(v)$ to induce a blue $C_4$.

*Proof of claim 1(b).* Suppose that Claim 1(b) is false. Assume that $v_{1,1}$ contains a red $(3,2)$ split. In particular, suppose that $v_{1,1}$ is adjacent in red to $\{v_{2,1},v_{2,2},v_{2,3},v_{3,1},v_{3,2}\}$.

Next, without loss of generality assume that $(v_{2,1},v_{3,2})$ is a red edge as $\{v_{2,1},v_{2,2},v_{3,1},v_{3,2}\}$ cannot induce a blue $C_4$. Next as $\{v_{2,2},v_{2,3},v_{3,1},v_{3,2}\}$ cannot induce a blue $C_4$ we may assume without loss of generality that $(v_{2,2},v_{3,2})$ is a red edge. Since $N_R(v_{1,1})$, cannot have a red $2K_2$, we may assume that $v_{3,1}$ is adjacent in blue to $v_{2,1},v_{2,2}$ and $v_{2,3}$.

This gives rise to two possible scenarios. In the first scenario, assume $v_{1,2}$ is adjacent in red to $v_{3,2}$. But then in order to avoid a red $B$, $(v_{1,2},v_{2,2})$ is a blue edge as illustrated in the following figure.

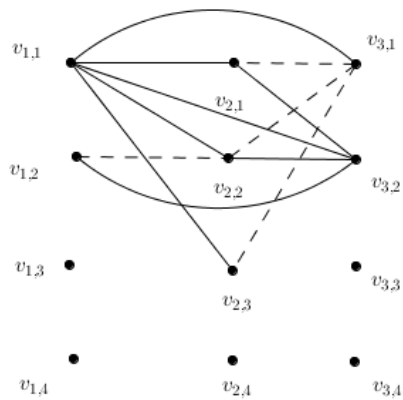

Figure 3.7: Graph used to show claim 1(b)-scenario 1

But in order to avoid a blue $C_4$, $(v_{1,2}, v_{2,1})$ will have to be a red edge. This will give us a red $B$, with centre $v_{3,2}$, a contradiction.

Next in the second scenario, as we have eliminated the first scenario, we may assume that $(v_{1,2}, v_{3,2})$, $(v_{1,3}, v_{3,2})$ and $(v_{1,4}, v_{3,2})$ are blue edges. Also in order to avoid a blue $C_4$, without loss of generality we may assume that at least two of $(v_{1,2}, v_{2,2})$, $(v_{1,3}, v_{2,2})$ and $(v_{1,4}, v_{2,2})$ edges are red edges. Thus, without loss of generality we may assume that $(v_{1,2}, v_{2,2})$ and $(v_{1,3}, v_{2,2})$ are red edges. Also in order avoid a blue $C_4$, $(v_{1,3}, v_{2,1})$ or $(v_{1,2}, v_{2,1})$ must be a red edge. Without loss of generality, in order to avoid a blue $C_4$, assume that $(v_{1,3}, v_{2,1})$ is a red edge. Next note that if $v_{3,3}$ is adjacent to both $v_{1,3}$ and $v_{1,4}$, we get a blue $C_4$. Therefore, without loss of generality we may assume that $(v_{1,3}, v_{3,3})$ is a red edge.

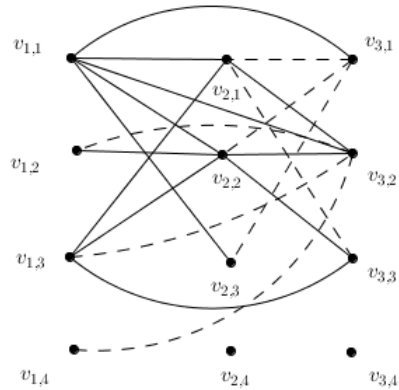

Figure 3.8: Graph used to show claim 1(b) scenario 2

In order to avoid a red $B$, $(v_{2,1}, v_{3,3})$ will have to be a blue edge. But to avoid a blue $C_4$, $(v_{3,3}, v_{2,2})$ will have to be a red edge as illustrated in the above figure. But this gives us a red $B$, with centre $v_{2,2}$, a contradiction.

Continuing with the main part of the proof of $j = 3$ case, applying claim 1(b), we get that $v_{1,1}$ give rise to a blue, (2,2) or contain a blue (3,0) split.

First suppose that, no vertex of $G$ has a blue (2,2) split. Then as each of the 4 vertices of $V_1$ are adjacent to at least 3 vertices of $V_2$ in blue or is adjacent to at least 3 vertices of $V_3$, in blue. By Pigeon Hole principle we will get that there are two vertices in $V_1$ having two common blue neighbours of $V_1^c$. This will force a blue $C_4$, a contradiction. Therefore, we may assume that a vertex of $V_1$ has a blue (2,2) split and in particular, we may assume that, $v_{1,1}$ is adjacent in blue to $U_1 = \{v_{2,1}, v_{2,2}, v_{3,1}, v_{3,2}\}$. This will give rise to two cases.

*Case 1:* If $v_{1,2}$, $v_{1,3}$ and $v_{1,4}$ contains (3,0) blue splits.

In particular we may assume that, $v_{1,2}$ is adjacent to $\{v_{2,2}, v_{2,3}, v_{2,4}\}$. However, as $v_{1,3}$, $v_{1,4}$ vertices also contain a (3,0) splits and they will force $v_{1,3}$ and $v_{1,4}$ to be adjacent to three vertices each of $V_3$. This will force a blue $C_4$, a contradiction.

*Case 2:* Without loss of generality if $v_{1,2}$ has a blue (2,2) split.

Without loss of generality we may assume that $v_{1,2}$ is adjacent to $\{v_{2,2}, v_{2,3}, v_{3,3}, v_{3,4}\}$ (if it is adjacent in blue to $\{v_{2,2}, v_{2,3}, v_{3,2}, v_{3,3}\}$ it will force a blue $C_4$ and if it is adjacent in blue to $\{v_{2,3}, v_{2,4}, v_{3,3}, v_{3,4}\}$, it will force $v_{1,3}$ to lie in a blue $C_4$).

In this situation, in order to avoid a blue $C_4$, $v_{1,3}$ (or $v_{1,4}$) cannot contain a blue (2,2) split. Therefore, both $v_{1,3}$ and $v_{1,4}$ must contain blue (3,0) splits. In this situation firstly $v_{1,3}$ and $v_{1,4}$ cannot be adjacent to $v_{2,2}$. Also both $v_{1,3}$ and $v_{1,4}$ can be adjacent to at most 2 vertices of $V_3$ in blue.

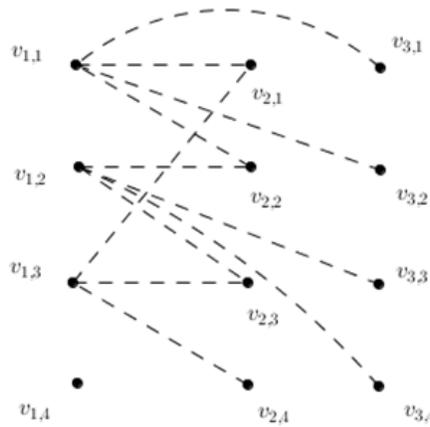

Figure 3.9: Blue graph generated by case 1

Therefore, as illustrated in the above figure $v_{1,3}$ must be adjacent to $v_{2,1}$, $v_{2,3}$ and $v_{2,4}$. Similarly, $v_{1,4}$ also must be adjacent to $v_{2,1}$, $v_{2,3}$ and $v_{2,4}$. But this will force a blue $C_4$, a contradiction. Therefore, we get that, $m_3(B,C_4) \leq 4$. That is, $m_3(B,C_4) = 4$ as required.

If $j = 4$, consider the coloring $G = K_{4\times 2} = H_R \oplus H_B$ of $K_{4\times 2}$ where $H_R$ and $H_B$ are the red and blue subgraph of $G$ induced by the red and blue colorings such that $H_R$ consists of a four-regular graph as illustrated in the following figure.

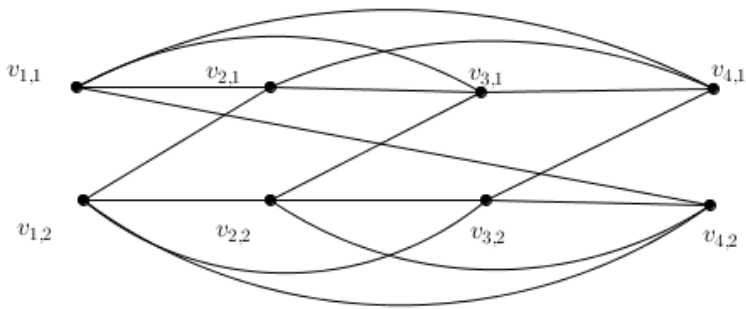

Figure 3.10: $H_R$ graph used to show $m_4(B,C_4) \geq 3$.

Clearly, $H_R$ has no red $B$ as any pair of red triangles are disjoint or else have an edge in common. $H_B$ is isomorphic to a blue $C_8$ and thus will not contain a blue $C_4$. Therefore, $m_4(B,C_4) \geq 3$. Next to show $m_4(B,C_4) \leq 3$, consider any red/blue coloring given by $G = K_{4\times 3} = H_R \oplus H_B$, such that $H_R$ contains no red $B$ and $H_B$ contains no blue $C_4$.

*Claim 2: a)* No vertex of $G$ can have red degree equal to 6 or more.

*b)* All vertices of $H_B$ will have blue degree 4.

*Proof of claim 2(a)* Assume the claim 2(a) is false. Without loss of generality we may assume that, $v_{1,1}$ is adjacent in red to exactly 6 vertices as if it is adjacent to 7 or more vertices the result will be clearly true. This situation will give rise to a red, (3,3,0), (3,2,1) or (2,2,2) split. In all these scenarios as $N_R(v_{1,1})$ cannot induce a red $2K_2$. Therefore, in all these scenarios, by inspection, $N_R(v_{1,1})$ will induce a blue $C_4$, a contradiction. Hence, claim 2(a) is true.

*Proof of claim 2(b)* Assume the claim 2(b) is false. Without loss of generality we may assume that, $v_{1,1}$ is adjacent in blue to 5 or more vertices. By claim 2(a), both $v_{1,2}$ and $v_{1,3}$ are adjacent in blue to 4 or more vertices. By Pigeon Hole principle as $|V_1^c| = 9$, we see that for this to happen at least two vertices of $V_1$ must have at least 2 common blue neighbours in $V_2 \cup V_3 \cup V_4$. This will result in a blue $C_4$, a contradiction. Hence, we get claim 2(b).

Continuing with the main part of the proof of $j = 4$ case, applying claim 2(a) and claim 2(b), we get that $v_{1,1}$ give rise to a blue, (2,2,0), (2,1,1) or (3,1,0) split. In the first scenario, without loss of generality, we may assume that $v_{1,1}$ is adjacent in blue to $U = \{v_{2,1},v_{2,2},v_{3,1},v_{3,2}\}$ and is adjacent in red to $W = \{v_{2,3},v_{3,3},v_{4,1},v_{4,2},v_{4,3}\}$. But then, as $W$ has no red $2K_2$, $W$ will be forced to have a blue $C_4$, a contradiction. In the second scenario, without loss of generality, we may assume that, $v_{1,1}$ has a red, (2,2,1) split and is adjacent in red to $U_1 = \{v_{2,1},v_{2,2},v_{3,1},v_{3,2},v_{4,1}\}$. But then, in order to avoid a red $B$, $U_1$ cannot contain a red $2K_2$ and this will force $U_1$ to induce a blue $C_4$ which would result in a contradiction, unless without loss of generality $V_1 = \{v_{2,2},v_{3,2},v_{4,1}\}$ induce a red $C_3$ and all the other 5 edges induced by $U_1$ are blue. In such a situation, as $v_{1,1}$ is adjacent in blue to $V_1^C \cap U_1^C$ and the blue degree of $v_{1,2}$ and $v_{1,3}$ are four, we get that both $v_{1,2}$ and $v_{1,3}$ are adjacent in blue to $v_{2,2}$ and $v_{3,2}$. This gives us a blue $C_4$, a contradiction.

Therefore, we see that vertex $v_{1,1}$ cannot give rise to a blue, (2,2,0) split or a (2,1,1) split. In fact, this argument can be extended to any arbitrary vertex of $G$. That is, no vertex of $G$ can have either a blue, (2,2,0) split or a (2,1,1) split. In other words, we can assume that all vertices of $G$ must have a red, (3,2,0) split (that is a blue (3,1,0) split). Therefore, without loss of generality we may assume that $v_{1,1}$ is adjacent to all vertices of $V_2$ in blue, $v_{1,2}$ is adjacent to all vertices of $V_3$ in blue and $v_{1,3}$ is adjacent to all vertices of $V_4$ in blue. However, as $v_{1,1}$ is adjacent to all but one vertex of $V_3 \cup V_4$ in red, there must be at least two adjacent blue edges among $V_3$ and $V_4$. Without loss of generality, assume that the two adjacent blue edges are given

# A RAMSEY PROBLEM RELATED TO BUTTERFLY GRAPH VS. PROPER CONNECTED SUBGRAPHS OF $K_4$

by $(x,y)$ and $(x,z)$ where $x \in V_3$ and $y,z \in V_4$. But then, $y\,x\,z\,v_{1,3}\,y$ will be a blue $C_4$, a contradiction. Therefore, we get $m_4(B,C_4) \leq 3$. That is, $m_4(B,C_4) = 3$.

In the $j \in \{5,6\}$ case, first consider the coloring $K_{6 \times 1} = H_R \oplus H_B$ of $K_{6 \times 1}$ where $H_R$ and $H_B$ are the red and blue subgraph of $G$ induced by the red and blue colorings such that $H_R$ consists of a $K_{3,3}$ and $H_B$ consists of a $2K_3$. Then the graph has no red $B$ and has no blue $C_4$. Therefore, $m_6(B,C_4) \geq 2$. Next to show $m_5(B,C_4) \leq 2$, consider any red/blue coloring given by $G = K_{5 \times 2} = H_R \oplus H_B$ such that $H_R$ contains no red $B$ and $H_B$ contains no blue $C_4$.

*Claim 3: a)* No vertex of $G$ can have red degree equal to 5 or more. In addition, $G$ contains no red $K_4$.

*b)* Let $T = \{x_1, x_2, y_1, y_2\}$ represent any four elements of $G = K_{5 \times 2}$ where $\{x_1, x_2\}$ belong to one partite set and $\{y_1, y_2\}$ belong to another partite set. Then there can be at most two red edges induced by $T$.

*c)* Let $T = \{x_1, x_2, y_1, y_2\}$ represent any four elements of $G = K_{5 \times 2}$ where $\{x_1, x_2\}$ belong to one partite set and $\{y_1, y_2\}$ belong to another partite set. If there are exactly two edges two red edges induced by $T$ then these two edges must be adjacent to each other.

*Proof of claim 3(a)* Assume the first part of the claim 3(a) is false. Without loss of generality we may assume that $v_{1,1}$ is adjacent in red to at least all vertices of $U = \{v_{2,1}, v_{2,2}, v_{3,1}, v_{4,1}, v_{5,1}\}$ or else adjacent in red to at least all vertices of $W = \{v_{2,1}, v_{2,2}, v_{3,1}, v_{3,2}, v_{4,1}\}$. In the first scenario as there is no red $2K_2$ clearly $U$ will induce a blue $C_4$ unless $U$ contains a red $C_3$. Without loss of generality, suppose that the red $C_3$ is given by $v_{2,1} v_{3,1} v_{4,1} v_{2,1}$. However, even in this case we will get a blue $C_4$ induced by $\{v_{2,2}, v_{3,1}, v_{4,1}, v_{5,1}\}$. Therefore, we are left with the second scenario. In the second scenario as there is no red $2K_2$ among $W = \{v_{2,1}, v_{2,2}, v_{3,1}, v_{3,2}, v_{4,1}\}$.

Thus, $W$ will induce a blue $C_4$ unless $W$ contains a red $C_3$. Without loss of generality, suppose that the red $C_3$ is given by $v_{2,1} v_{3,1} v_{4,1} v_{2,1}$. This would force $\{v_{2,2}, v_{3,2}, v_{4,1}\}$ to induce a blue $C_3$. As $\{v_{2,1}, v_{3,1}, v_{5,1}, v_{5,2}\}$ will not induce a blue $C_4$, without loss of generality we may assume that $(v_{5,1}, v_{3,1})$ is a red edge. Next in order to avoid a red $B$ with centre $v_{3,1}$, $(v_{5,1}, v_{4,1})$ is a blue edge. Also $(v_{2,1}, v_{3,2})$ is a blue edge in order to avoid a red $B$ with centre $v_{1,1}$. Next in order to avoid a blue $C_4$, $(v_{2,1}, v_{5,1})$ is a red edge. This results in a red $B$ with centre $v_{3,1}$ (consisting of the two red triangles $v_{2,1} v_{3,1} v_{5,1} v_{2,1}$ and $v_{1,1} v_{3,1} v_{4,1} v_{1,1}$). Hence, the first part of claim 3(a) follows. Next assume that there is a red $K_4$ say induced by $X = \{v_{1,1}, v_{2,1}, v_{3,1}, v_{4,1}\}$. Then in order to avoid a red $B$ both $v_{5,1}$ and $v_{5,2}$ will be adjacent at most 1 vertex of $X$ in red. But then this will force both $v_{5,1}$ and $v_{5,2}$ to be adjacent in blue to three vertices each of $T$. This will give us a blue $C_4$ containing $v_{5,1}$ and $v_{5,2}$, a contradiction.

*Proof of claim 3(b)* Next assume the claim 3(b) is false. Let $T = \{x_1, x_2, y_1, y_2\}$. Suppose that there are at least 3 red edges induced by $T$. Then this will result in a red $P_4$ in $T$, such that two vertices have red degree 2 in $T$. Denote these two vertices by $x$ and $y$. But then by part 3(a), both $x$ and $y$ will have two common blue neighbours in $T^c$. This will force a blue $C_4$, a contradiction.

*Proof of claim 3(c)* Next assume the claim 3(c) is false. Without loss of generality, assume that both the red induced subgraph of $V_1 \cup V_2$ consists of exactly two red edges namely $(v_{1,1}, v_{2,2})$ and $(v_{1,2}, v_{2,1})$. Then as $m_4(C_3, C_4) = 2$ (see [6]) we get that there is a $C_3$ and it will give rise to two cases.

**Case 1:** The red $C_3$ is induced by $\{v_{2,1}, v_{3,1}, v_{4,1}\}$

Next in order to avoid a red $K_4$ induced by $\{v_{1,2}, v_{2,1}, v_{3,1}, v_{4,1}\}$, without loss of generality we may assume that $(v_{1,2}, v_{3,1})$ is a red edge and $(v_{1,2}, v_{4,1})$ is blue or else $(v_{1,2}, v_{3,1})$ and $(v_{1,2}, v_{4,1})$ are both blue edges. In the first scenario in order to avoid a blue $C_4$ induced by $\{v_{3,1}, v_{4,1}, v_{5,1}, v_{5,2}\}$, without loss of generality we may assume that $(v_{3,1}, v_{5,1})$ is a red edge. Next, in order to avoid a red $B$, $(v_{4,1}, v_{5,1})$ and $(v_{5,1}, v_{1,2})$ have to be blue edges. In order to avoid a blue $C_4$ induced by $\{v_{2,1}, v_{3,1}, v_{4,2}, v_{5,2}\}$ and as red degree of $v_{3,1}$ is four, without loss of generality we may assume that $(v_{2,1}, v_{5,2})$ is a red edge. Next, as red degree of $v_{2,1}$ and $v_{3,1}$ is also four, $(v_{2,1}, v_{4,2})$, $(v_{3,1}, v_{1,1})$, $(v_{3,1}, v_{4,2})$ have to be blue edges. This gives rise to the following figure.

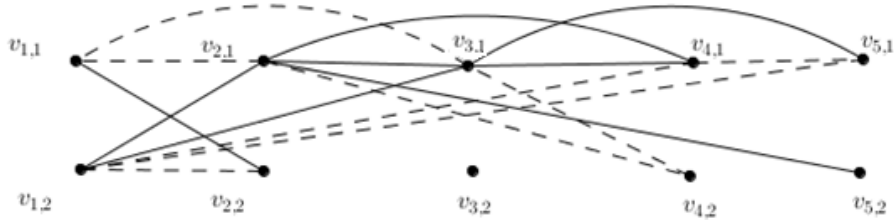

Figure 3.11: Graph used in proof of case 1-first scenareo

But then as seen in the figure, we will get that $\{v_{1,1}, v_{2,1}, v_{3,1}, v_{4,2}\}$ will induce a blue $C_4$, a contradiction.

In the second scenario, as $(v_{1,2}, v_{3,1})$ and $(v_{1,2}, v_{4,1})$ are both blue edges. In order to avoid a blue $C_4$ induced by $\{v_{1,1}, v_{1,2}, v_{3,1}, v_{4,1}\}$ without loss of generality we may assume that $(v_{1,1}, v_{3,1})$ is a red edge. Next in order to avoid a red $B$, $(v_{2,2}, v_{3,1})$ will have to be a blue edge. In order to avoid a blue $C_4$ we get that $(v_{2,2}, v_{4,1})$ is a red edge, as illustrated in the following figure.

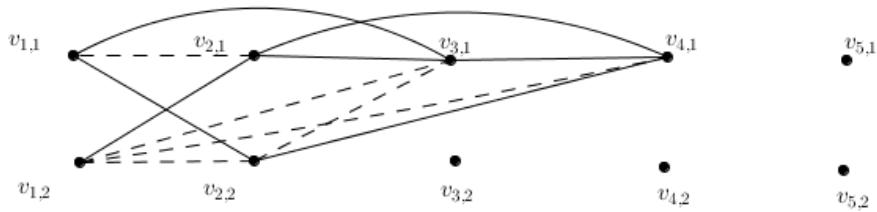

Figure 3.12: Graph used in proof of case 1-second scenario

As red degrees of $v_{2,1,}$ and $v_{4,1}$ are at most four, there, will be 2 vertices in $\{v_{1,1}, v_{3,2}, v_{5,1}, v_{5,2}\}$ which are adjacent to $v_{2,1,}$ and $v_{4,1}$ in blue. Suppose these two vertices are $x$ and $y$. Then $(y, v_{2,1}, x, v_{4,1}, y)$ will be a blue $C_4$, a contradiction.

# A RAMSEY PROBLEM RELATED TO BUTTERFLY GRAPH VS. PROPER CONNECTED SUBGRAPHS OF $K_4$

**Case 2:** The red $C_3$ is induced by $T = \{v_{3,1}, v_{4,1}, v_{5,1}\}$

In order to avoid case 1, we may assume that, each of the four vertices will be adjacent in red to at most one vertex of $T$. That is each of the four vertices of the set $\{v_{1,1}, v_{1,2}, v_{2,1}, v_{2,2}\}$ will be adjacent in blue to two vertices of the set $T$. This will force a pair of vertices of $\{v_{1,1}, v_{1,2}, v_{2,1}, v_{2,2}\}$ to have two common neighbours in $T$. Thus, we will get a blue $C_4$, a contradiction.

Continuing with the main proof of $j \in \{5,6\}$ case, from [6] we get that there is a red $C_3$ in $H_R$. Without loss of generality, assume that the red $C_3$, is induced by say $\{v_{1,1}, v_{2,1}, v_{3,1}\}$. Let $S = \{v_{1,1}, v_{2,1}, v_{3,1}\}$. Applying claim 3(b) and claim 3(c), we get that $S_1 = \{v_{1,2}, v_{2,2}, v_{3,2}\}$ will induce a blue $C_3$ as illustrated in the following diagram.

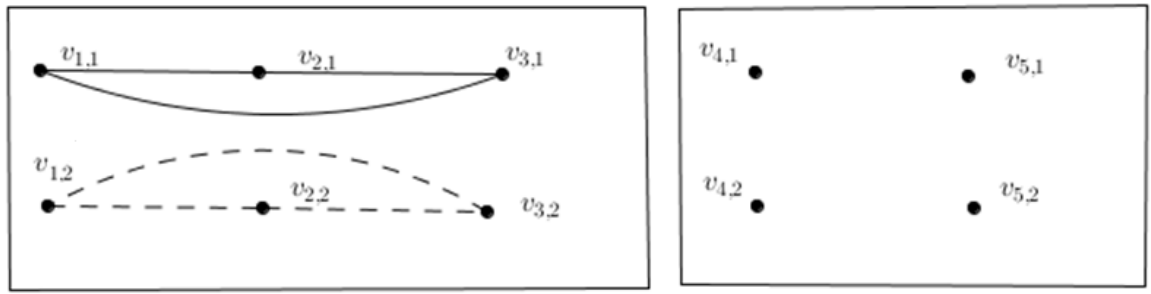

Figure 3.13: Graph used in proof of case 2.

Note that in order to avoid a blue $C_4$, either $(v_{2,1}, v_{3,2})$ or $(v_{2,1}, v_{1,2})$ must be red. Without loss of generality, assume that $(v_{2,1}, v_{3,2})$ is red. Then by claim 3(b), $(v_{2,2}, v_{3,1})$ must be blue. In order to avoid a blue $C_4$, this will force $(v_{1,2}, v_{3,1})$ to be red. That is, we get that both $v_{2,1}$ and $v_{3,1}$ are adjacent to at least three vertices of $S \cup S_1$ in red. Hence, $v_{2,1}$ and $v_{3,1}$ will have two common neighbour in blue in $\{v_{4,1}, v_{4,2}, v_{5,1}, v_{5,2}\}$. This results in a blue $C_4$, a contradiction. Therefore, we get $m_5(B, C_4) \leq 2$. That is, $m_6(B, C_4) = m_5(B, C_4) = 2$.

If $j \geq 7$, since $r(B, C_4) = 7$ (see Henery (1989)), we get $m_j(B, C_4) = 1$.